# CHOICE OF NEIGHBOR ORDER IN NEAREST-NEIGHBOR CLASSIFICATION

By Peter Hall, Byeong U. Park [1] and Richard J. Samworth

*Australian National University, Seoul National University and University of Cambridge*

The $k$th-nearest neighbor rule is arguably the simplest and most intuitively appealing nonparametric classification procedure. However, application of this method is inhibited by lack of knowledge about its properties, in particular, about the manner in which it is influenced by the value of $k$; and by the absence of techniques for empirical choice of $k$. In the present paper we detail the way in which the value of $k$ determines the misclassification error. We consider two models, Poisson and Binomial, for the training samples. Under the first model, data are recorded in a Poisson stream and are "assigned" to one or other of the two populations in accordance with the prior probabilities. In particular, the total number of data in both training samples is a Poisson-distributed random variable. Under the Binomial model, however, the total number of data in the training samples is fixed, although again each data value is assigned in a random way. Although the values of risk and regret associated with the Poisson and Binomial models are different, they are asymptotically equivalent to first order, and also to the risks associated with kernel-based classifiers that are tailored to the case of two derivatives. These properties motivate new methods for choosing the value of $k$.

**1. Introduction.** In the classification or discrimination problem with two populations, denoted by $X$ and $Y$, one wishes to classify an observation $z$ to either $X$ or $Y$ using only training data. The $k$th-nearest neighbor classification rule is arguably the simplest and most intuitively appealing nonparametric classifier. It assigns $z$ to population $X$ if at least $\frac{1}{2}k$ of the $k$ values in the pooled training-data set nearest to $z$ are from $X$, and to population $Y$ otherwise. The first study of this method was undertaken by Fix and Hodges

---

Received June 2007; revised August 2007.

[1]Supported by KOSEF (project R02-2004-000-10040-0).

(1951). Since then there have been many investigations into the method's statistical properties. Little is known about the structure of its error probabilities, however, and neither are formulae available for optimal choice of $k$. Practical methods for optimal empirical choice of $k$ have apparently not been given.

The present paper resolves these issues, and focuses on expansions of the error rate of $k$th-nearest neighbor classifiers which are associated with optimal choice of $k$. We show that the values of risk of nearest-neighbor classifiers can be represented quite simply in terms of properties of the two populations, and that this leads to new, practical ways of choosing the value of $k$.

The sizes of the training samples used to construct classifiers might fairly be viewed as random variables. Consider, for example, the case where a classifier is used by a bank to determine, from the bank's data, whether a new customer is likely to default on a loan. The sizes of the two training samples could be the number, $M$, of previous customers who defaulted, and the number, $N$, who did not default, respectively. An appropriate model for the distributions of $M$ and $N$ might be that they are statistically independent and Poisson, with means $\mu$ and $\nu$, say. For example, the Poisson sample-size model could arise if the population of potential customers were much larger than the number of customers who sought loans from the bank.

Thus, Poisson rather than deterministic models for training-sample sizes can be motivated. Here, the total number of data in the two training samples is random, and data in a Poisson stream are "assigned" to one or other of the two populations using a formula which is based on the respective prior probabilities. A different approach, which gives rise to a Binomial-type model, involves the total number of training data being pre-determined, but apportions these data among the two populations in a manner similar to the Poisson model. We shall show that these two approaches produce nearest-neighbor classifiers with risks that are different but are nevertheless first-order equivalent.

For fixed $k$ the risk of a $k$-nearest neighbor classifier converges to its limit relatively quickly, at rate $T^{-2}$, as total sample size, $T$, increases [Cover (1968)]. However, the limiting value is strictly larger than the Bayes risk of the "ideal" classifier that would be used if both population densities were known. By way of comparison, in the case of imperfect information about the population, and in particular, in parametric settings, the risk of empirical Bayes classifiers converges to the Bayes risk no more rapidly than $T^{-1}$; see Kharin and Duchinskas (1979). In nonparametric settings the rate of convergence to Bayes risk is slower still, but may nevertheless be asymptotically optimal; see, for example, Marron (1983) and Mammen and Tsybakov (1999).

In most previous work on nearest-neighbor classifiers, the value of $k$ was held fixed. Cover and Hart (1967) gave upper bounds for the limit of the risk of nearest-neighbor classifiers. Wagner (1971) and Fritz (1975) treated convergence of the conditional error rate when $k = 1$. Devroye and Wagner (1977, 1982) developed and discussed theoretical properties, particularly issues of mathematical consistency, for $k$-nearest-neighbor rules. Devroye (1981) found an asymptotic bound for the regret with respect to the Bayes classifier. Devroye et al. (1994) gave a particularly general description of strong consistency for nearest-neighbor methods. Psaltis, Snapp and Venkatesh (1994) generalized the results of Cover (1968) to general dimension, and Snapp and Venkatesh (1998) further extended the results to the case of multiple classes. Bax (2000) gave probabilistic bounds for the conditional error rate in the case where $k = 1$. Kulkarni and Posner (1995) addressed nearest-neighbor methods for quite general dependent data, and Holst and Irle (2001) provided formulae for the limit of the error rate in the case of dependent data. Related research includes that of Györfi (1978, 1981) and Györfi and Györfi (1978), who investigated the rate of convergence to the Bayes risk when $k$ tends to infinity as $T$ increases.

In the case of classifiers based on second-order kernel density estimators, and for populations with twice-differentiable densities, the risk typically converges to the Bayes risk at rate $n^{-4/(d+4)}$, where $d$ denotes the number of dimensions. See, for example, Kharin (1982), Raudys and Young (2004) and Hall and Kang (2005). In a minimax sense that Marron (1983) makes precise, this rate is optimal. As we show in this paper, nearest-neighbor classifiers with Poisson or Binomial interpretations of sample size have the same property.

Recent work on properties of classifiers focuses largely on deriving upper and lower bounds to regret in cases where the classification problem is relatively difficult, for example, where the classification boundary is comparatively unsmooth. Research of Audibert and Tsybakov (2005) and Kohler and Krzyzak (2006), for example, is in this category. The work of Mammen and Tsybakov (1999), which permits the smoothness of a classification problem to be varied in the continuum, forms something of a bridge between the smooth case, which we treat, and the rough case.

There is a literature on empirical choice of $k$; see, for example, Chapter 26 of Devroye, Györfi and Lugosi (1996) and Sections 7.2 and 8.4 of Györfi et al. (2002). More generally, Devroye, Györfi and Lugosi (1996) explored the properties and features of nearest-neighbor methods in the setting of pattern recognition. Chapter 5 of that monograph gives a good guide to the literature in this setting.

## 2. Main results.

2.1. *Different interpretations of sample size.* Assume we have $m$ identically distributed data $\mathcal{X} = \{X_1, \ldots, X_m\}$, and $n$ identically distributed data $\mathcal{Y} = \{Y_1, \ldots, Y_n\}$, all of them $d$-variate and mutually independent. Let the respective probability densities be $f$ and $g$. Given a compact set $\mathcal{R} \subseteq \mathbb{R}^d$, we wish to use the data to classify a new datum $z \in \mathcal{R}$ as coming from the $X$ or $Y$ population. Note that we do not assume $f$ and $g$ themselves to be compactly supported; the constraint is only that we confine attention to the problem of classifying new data that come from a given compact region $\mathcal{R}$.

In many instances the ratio of the sizes of the datasets is a good approximation to the ratio of the prior probabilities of observing the respective populations. We shall adopt this viewpoint, which raises the issue of how we should interpret $m$ and $n$. Two models arise in a natural way: the Poisson, where the individual sample sizes are Poisson-distributed and data are assigned randomly to one proportion or another, in proportion to the respective likelihoods; and the Binomial, where the sum of the two training-sample sizes is deterministic but data are ascribed to populations in the same fashion as before. The Poisson case can be viewed as the result of taking a sample from a marked point process in $\mathbb{R}^d$, and assigning marks in a way that reflects prior probabilities; and the Binomial case is the result of conditioning on total sample size in the Poisson setting.

In the sense that it avoids the conditioning step, the Poisson case is the more natural and has the greater degree of symmetry. Therefore, we take that as the basis for analysis, and tackle the Binomial model by reference to the solution in the Poisson case.

In multi-population cases, the $k$th nearest-neighbor classifier would typically be used to assign $z$ to population $j$ if that population accounted for the greatest number of data among the $k$ values in the pooled dataset that are nearest to $z$. Our results apply directly to this case, provided we work within a compact region at each point of which the maximum value of the population densities is achieved by no more than two densities. Another straightforward extension is to the case where distance is measured in a weighted Euclidean metric; we shall work only with the standard, unweighted form.

2.2. *Poisson model.* Assume that $\mathcal{X} = \{X_1, X_2, \ldots\}$ and $\mathcal{Y} = \{Y_1, Y_2, \ldots\}$ represent points of type $X$ and type $Y$, respectively, in a two-type marked Poisson process, $\mathcal{P}$, in $\mathbb{R}^d$, with intensity function $\mu f + \nu g$, and respective probabilities

$$(2.1) \qquad\qquad \psi(z) = \frac{\mu f(z)}{\mu f(z) + \nu g(z)}$$

and $1 - \psi(z)$ that a point of $\mathcal{P}$ at $z$ is of type $X$ or of type $Y$. In particular, the respective prior probabilities of the $X$ and $Y$ populations are $\mu/(\mu + \nu)$ and $\nu/(\mu + \nu)$. It will be assumed that $f$ and $g$ are held fixed, and that $\mu$ and $\nu$ satisfy

(2.2)    $\mu = \mu(\nu)$ increases with $\nu$, in such a manner that $\mu/(\mu + \nu) \to p \in (0, 1)$ as $\nu \to \infty$.

Define $\rho = pf/\{pf + (1 - p)g\}$, a function on $\mathbb{R}^d$.

Suppose too that the respective densities, $f$ and $g$, of the $X$ and $Y$ populations, satisfy

(2.3)    the set $\mathcal{S} \subseteq \mathcal{R}$, defined as the locus of points $z$ for which $\rho(z) = \frac{1}{2}$, is of codimension 1 and of finite measure in $d - 1$ dimensions.

(2.4)    the distributions with densities $f$ and $g$ have finite second moments; $f$ and $g$ are both continuous in an open set containing $\mathcal{R}$, and both have two continuous derivatives within an open set containing $\mathcal{S}$; and $f + g > 0$ on $\mathcal{R}$;

The first part of (2.3) asks that $\mathcal{S}$ be a $(d - 1)$-dimensional structure—a set of isolated points if $d = 1$, a set of curves in the plane if $d = 2$, and so on.

The assumption of two derivatives in (2.4) is to be expected since, as noted in Section 1, the convergence rate of regret that is achieved by nearest-neighbor methods is optimal under that smoothness assumption. The condition that the derivatives assumed in (2.4) are continuous is imposed only so that a concise asymptotic formula for regret can be given; see (2.8) below. Without the precision provided by the continuity assumption, we could state only an upper bound for regret, in which the right-hand side of (2.8) was replaced by $O\{k^{-1} + (k/\nu)^{4/d}\}$.

We ask too that the slopes at which the two densities, weighted in proportion to their prior probabilities, meet along $\mathcal{S}$ be bounded away from zero along $\mathcal{S}$. That is, the function

$$a(z)^2 \equiv \sum_{j=1}^{d} \left\{ p \frac{\partial f(z)}{\partial z_j} - (1 - p) \frac{\partial g(z)}{\partial z_j} \right\}^2 \text{ is bounded away from zero on } \mathcal{S}.$$

(2.5)

Equivalently, the prior-weighted densities cross at an angle, rather than meet in a tangential way. If the prior-weighted densities were to have exactly equal gradients at crossing places, then there would be an explicit and intimate connection between the distributions of $X$ and $Y$ populations that could hardly arise by chance. It is difficult to envisage that perfect alignment of densities at crossing points would actually occur commonly in practice.

Write $dz_0$ for an infinitesimal element of $\mathcal{S}$, centred at $z_0$. Let $a_d = \pi^{d/2}/\Gamma(1 + \frac{1}{2}d)$ denote the content of the unit $d$-dimensional sphere, define

$\lambda = p(1-p)^{-1}f + g$ and

$$(2.6) \quad \alpha(z) = \frac{d}{d+2}\lambda(z)^{-1-(2/d)}a_d^{-2/d}d^{-1}\sum_{j=1}^d \left\{\rho_j(z)\lambda_j(z) + \frac{1}{2}\rho_{jj}(z)\lambda(z)\right\},$$

where $z = (z^{(1)}, \ldots, z^{(d)})$, $\lambda_j(z)$ and $\rho_j(z)$ denote the first derivatives of the respective functions with respect to $z^{(j)}$, and $\rho_{jj}(z)$ is the second derivative of $\rho(z)$ with respect to $z^{(j)}$. Put $\dot\rho = (\rho_1, \ldots, \rho_d)$.

Let $\Phi$ denote the cumulative distribution function of the standard normal distribution and let $\Psi_1(z) = \|\dot\rho(z)\|^{-1}\{p\dot f(z) - (1-p)\dot g(z)\}^{\mathrm{T}}\dot\rho(z)$. It can be shown that, on $\mathcal{S}$, $a(z) = \Psi_1(z) = 4h(z)\|\dot\rho(z)\|$, where $h(z)$ denotes the common value that $pf(z)$ and $(1-p)g(z)$ assume at $z \in \mathcal{S}$. Therefore, since assumptions (2.3)–(2.5) imply that $a$ and $h$ are bounded away from zero and infinity on $\mathcal{S}$, they also ensure that $\Psi_1(z)$ and $\|\dot\rho(z)\|$ are bounded away from zero and infinity there. It follows that the constants $C_1$ and $C_2$, given by

$$(2.7) \quad \begin{aligned} C_1 &= \int_{\mathcal{S}} \Psi_1(z_0)(8\|\dot\rho(z_0)\|^2)^{-1}\,dz_0 = \frac{1}{2}\int_{\mathcal{S}}\frac{h(z_0)}{\|\dot\rho(z_0)\|}\,dz_0, \\ C_2 &= 4\int_{\mathcal{S}} \Psi_1(z_0)(8\|\dot\rho(z_0)\|^2)^{-1}\alpha(z_0)^2\,dz_0 = 2\int_{\mathcal{S}}\frac{h(z_0)}{\|\dot\rho(z_0)\|}\alpha(z_0)^2\,dz_0, \end{aligned}$$

are finite, that $C_1$ is nonzero, and that $C_2 = 0$ if and only if $\alpha$ is identically zero on $\mathcal{S}$.

The Bayes classifier assigns $z$ to the $X$ or $Y$ population according as $\psi(z) \geq \frac{1}{2}$ or $\psi(z) < \frac{1}{2}$, respectively. Therefore, the Bayes risk for classification on $\mathcal{R}$ is

$$\mathrm{risk}_{\mathrm{Bayes}}^{\mathrm{Pois}} = \int_{\mathcal{R}}\min\left(\frac{\mu f}{\mu+\nu}, \frac{\nu g}{\mu+\nu}\right),$$

where, here and below, the superscript "Pois" will indicate that the setting of the Poisson model is being considered. The risk of the $k$-nearest neighbor classifier, which assigns $z$ to population $X$ if at least $\frac{1}{2}k$ of the $k$ values of Poisson data nearest to $z$ are from $X$, and to population $Y$ otherwise, is

$$\mathrm{risk}_{k\text{-nn}}^{\mathrm{Pois}} = \frac{\mu}{\mu+\nu}\int_{\mathcal{R}} f(z)P^{\mathrm{Pois}}(z \text{ classified by } k\text{-nn rule as type } Y)\,dz$$
$$\qquad + \frac{\nu}{\mu+\nu}\int_{\mathcal{R}} g(z)P^{\mathrm{Pois}}(z \text{ classified by } k\text{-nn rule as type } X)\,dz.$$

A proof of the following result is given in Section 4.

THEOREM 1. *Assume the Poisson model, that (2.2)–(2.5) hold, and that $1 \leq k_1(\nu) < k_2(\nu)$, where $k_1(\nu)/\nu^\varepsilon \to \infty$ and $k_2(\nu) = O(\nu^{1-\varepsilon})$ for some $0 <$*

$\varepsilon < 1$. *Then,*

$$(2.8) \quad \text{risk}_{k\text{-nn}}^{\text{Pois}} - \text{risk}_{\text{Bayes}}^{\text{Pois}} = C_1 k^{-1} + C_2 (k/\nu)^{4/d} + o\{k^{-1} + (k/\nu)^{4/d}\},$$

*uniformly in* $k_1(\nu) \le k \le k_2(\nu)$.

Result (2.8) implies that, provided $\alpha$ is not identically zero, the optimal $k$ satisfies $k_{\text{opt}}^{\text{Pois}} \sim \text{const.} \nu^{4/(d+4)}$. To set (2.8) into context, we note that a general formula for the difference between the risk of an empirical classifier and the Bayes risk can be developed from the theory of "plug-in decisions"; see Theorem 2.2, page 16, of Devroye, Györfi and Lugosi (1996), and Theorem 6.2, page 93, of Györfi et al. (2002). When specialized to the case of nearest-neighbor methods, this argument bounds the left-hand side of (2.8) by a constant multiple of $\{k^{-1} + (k/\nu)^{2/d}\}^{1/2}$, the minimum order of which is $\nu^{-1/(d+2)}$. Mammen and Tsybakov (1999) showed that, in the case where discrimination boundaries are smooth, substantially faster convergence rates are possible. Result (2.8) and its analogues in the setting of Theorem 2 give concise accounts of those faster rates in the case of nearest-neighbor methods.

Expansion (2.8) has a close analogue in the setting of second-order, kernel-based methods. See, for example, formulae (3) of Kharin (1982) and (A.2) of Hall and Kang (2005).

2.3. *Binomial model.* In the Poisson model we can think of the data as arriving in a stream $(Z_1, L_1), (Z_2, L_2), \ldots$, where $Z_1, Z_2, \ldots$ comprise a Poisson process in $\mathbb{R}^d$, with intensity function $\mu f + \nu g$, and the "labels" $L_i$ form a sequence of zeros and ones, independent of one another conditional on the $Z_i$'s, with $P(L_i = 0 \mid Z_i) = \psi(Z_i)$ and $\psi$ defined by (2.1). If $L_i = 0$, then $Z_i$ is labeled as coming from the $X$-population, whereas if $L_i = 1$, then $Z_i$ is labeled as $Y$. Since the integral of the Poisson-process intensity over $\mathbb{R}^d$ equals $\mu + \nu$, then the number of points $Z_i$ equals a Poisson-distributed random variable, $T$, say, with mean $\mu + \nu$. In the Binomial model we use the same process to generate data, but now we condition on $T$.

It is convenient to think of $T$ as $m + n$, where $m = \mu T/(\mu + \nu)$ and $n = \nu T/(\mu + \nu)$ are the respective average numbers of points that would occur in the two training samples if we were to adopt the procedure indicated above. (In particular, $m$ and $n$ are not necessarily integers.) In this notation the risk for the nearest-neighbor classifier under the Binomial model can be written as

$$\text{risk}_{k\text{-nn}}^{\text{Bin}} = \frac{m}{T} \int_{\mathcal{R}} f(z) P^{\text{Bin}}(z \text{ classified by } k\text{-nn rule as type } Y) \, dz$$
$$+ \frac{n}{T} \int_{\mathcal{R}} g(z) P^{\text{Bin}}(z \text{ classified by } k\text{-nn rule as type } X) \, dz,$$

where we use the superscript Bin to indicate that we are sampling under the Binomial model. If we suppose that

$$(2.9) \qquad \mu + \nu = T \qquad \text{a nonrandom integer,}$$

then these manipulations are unnecessary, and so we shall assume (2.9) below. This condition also implies that the Bayes risk under the Binomial model, $\text{risk}_{\text{Bayes}}^{\text{Bin}}$, is identical to its counterpart under the Poisson model, and that helps to further simplify comparisons.

THEOREM 2. *Assume the Binomial model, that* (2.2)–(2.5) *and* (2.9) *hold, and that* $k_1$ *and* $k_2$ *satisfy the conditions imposed on them in Theorem* 1. *Then,*

$$(2.10) \qquad \text{risk}_{k\text{-nn}}^{\text{Bin}} - \text{risk}_{k\text{-nn}}^{\text{Pois}} = o\{k^{-1} + (k/\nu)^{4/d}\},$$

*uniformly in* $k_1(\nu) \le k \le k_2(\nu)$.

A proof of Theorem 2 is given in a longer version of this paper [Hall, Park and Samworth (2007)].

Formula (2.10) asserts that the difference between $\text{risk}_{k\text{-nn}}^{\text{Bin}}$ and $\text{risk}_{k\text{-nn}}^{\text{Pois}}$ is of smaller order than the difference between $\text{risk}_{k\text{-nn}}^{\text{Pois}}$ and $\text{risk}_{\text{Bayes}}^{\text{Pois}}$ [see (2.9) for the latter difference], and hence, implies that the expansion of regret at (2.8) is equally valid if $\text{risk}_{k\text{-nn}}^{\text{Pois}}$ and $\text{risk}_{\text{Bayes}}^{\text{Pois}}$ there are replaced by $\text{risk}_{k\text{-nn}}^{\text{Bin}}$ and $\text{risk}_{\text{Bayes}}^{\text{Bin}}$, respectively.

2.4. *Empirical choice of* $k_{\text{opt}}$. The theoretical results described earlier can be used to motivate practical methods for choosing $k$. We shall treat the Poisson model; the Binomial model can be addressed similarly.

Let $M$ and $N$ be the respective sizes of the training samples $\mathcal{X}$ and $\mathcal{Y}$. Generate $M^*$ and $N^*$, respectively, from the Poisson distributions with means equal to $M$ and $N$. Let $0 < r < 1$. Draw bootstrap resamples $\mathcal{X}^*$ and $\mathcal{Y}^*$, of respective sizes

$$M_1^* = [rM^*], \qquad N_1^* = [rN^*]$$

from $\mathcal{X}$ and $\mathcal{Y}$. Here, $[x]$ denotes the integer part of $x$. This choice of $M_1^*$ and $N_1^*$ implies that the total resample size equals $r(M^* + N^*)$, except for rounding errors arising from taking integer parts. Note too that $M_1^*/(M_1^* + N_1^*)$ equals the sampling fraction $M^*/(M^* + N_1^*)$ (again modulo integer-part rounding). This is necessary if our bootstrap algorithm, based on repeated resamples of sizes $M_1^*$ and $N_1^*$, is to mimic properties of the original sampling algorithm.

Draw additional resamples $\mathcal{X}_{\text{test}}^*$ and $\mathcal{Y}_{\text{test}}^*$, of respective sizes $M^* - M_1^*$ and $N^* - N_1^*$ from $\mathcal{X}$ and $\mathcal{Y}$. Build near-neighbor classifiers based on $\mathcal{X}^*$ and $\mathcal{Y}^*$. Use them to classify the data $\mathcal{X}^*_{\text{test}}$ and $\mathcal{Y}^*_{\text{test}}$, and compute the resulting error rate. Average this rate over a large number of choices of $\{\mathcal{X}^*, \mathcal{X}^*_{\text{test}}\}$ and $\{\mathcal{Y}^*, \mathcal{Y}^*_{\text{test}}\}$. Choose $k = \hat{k}_{\text{opt}}$ to minimize the average error rate; it is an estimator of the value of $\hat{k}_{\text{opt}}(r\mu, r\nu)$ that we would use if the true intensity function were $r(\mu f + \nu g)$, rather than $\mu f + \nu g$. Convert $\hat{k}_{\text{opt}}$ to an empirical value, $\tilde{k}_{\text{opt}} = r^{-4/(d+4)}\hat{k}_{\text{opt}}$, that is of the right size for classification starting from the samples $\mathcal{X}$ and $\mathcal{Y}$.

In the case of the binomial sample-size model, one may follow the same bootstrapping procedure as in the Poisson case, but generating $M^*$ from Binomial$(M + N, M/(M + N))$ and taking $N^* = M + N - M^*$.

**3. Numerical properties.** We present the results of a numerical experiment demonstrating the effectiveness of the empirical choice $\tilde{k}_{\text{opt}}$ introduced in Section 2. We simulated 500 training datasets from Poisson sample-size models for selected pairs of intensity constants $(\mu, \nu)$. Each dataset was obtained as follows. First, we generated a random number, say, $\mathcal{N}$, from a Poisson distribution with mean $\mu + \nu$. Then, we drew $\mathcal{N}$ independent data from the density $\lambda(z) = \{\mu f(z) + \nu g(z)\}/(\mu + \nu)$; let these be $Z_1, \ldots, Z_{\mathcal{N}}$. For $i = 1, \ldots, \mathcal{N}$, we marked "type $X$" or "type $Y$" on $Z_i$ with respective probabilities $\psi(Z_i) = \mu f(Z_i)/\{\mu f(Z_i) + \nu g(Z_i)\}$ and $1 - \psi(Z_i)$. An equivalent way of doing this would be to draw $\mathcal{N}$ independent data, each of which is sampled from the density $f$ or $g$ with respective probabilities $\mu/(\mu + \nu)$ and $\nu/(\mu + \nu)$. Each datum would then be marked "type $X$" if it was from $f$, and "type $Y$" otherwise.

We took $(\mu, \nu) = (100, 100)$ and $(100, 200)$ and considered the cases $d = 1, 2$. For $d = 1$, we chose $f$ to be the density function of $N(-0.5, 1)$ and $g$ to be the density function of $N(0.5, 1)$. For $d = 2$, we considered two pairs of densities. One was $(f, g)$, where $f \sim N_2((0.5, -0.5), I_2)$ and $g \sim N_2((-0.5, 0.5), I_2)$. Here, $I_d$ is the $d \times d$ identity matrix. The other was a pair of bivariate normal densities, as in the first case but with correlation $\rho = 0.5$.

For each $z$, we evaluated

$$\hat{P}^{\text{Pois}}(z \text{ classified as type } X)$$

$$= \tfrac{1}{500}(\# \text{ training samples that classify } z \text{ as type } X).$$

The error rate was then estimated by the formula

$$\widehat{\text{Err}} = \frac{\mu}{\mu + \nu} \int_{\mathcal{R}} f(z)\{1 - \hat{P}^{\text{Pois}}(z \text{ classified as type } X)\}\, dz$$

$$+ \frac{\nu}{\mu + \nu} \int_{\mathcal{R}} g(z)\hat{P}^{\text{Pois}}(z \text{ classified as type } X)\, dz.$$

We took $\mathcal{R} = [-2.5, 2.5]^d$, which covered most of the sampling region. To see the effect of the bootstrap resampling fraction on the performance of $\tilde{k}_{\text{opt}}$, the three choices $r = 1/3, 1/2, 2/3$ were considered, where $r$ was defined in Section 2.4. For computation of $\tilde{k}_{opt}$, 100 bootstrap resamples were drawn.

Table 1 shows the estimated error rates of the $k$-nearest neighbor classifier with $k_{opt}$, and the $k$-nearest neighbor classifier with $\tilde{k}_{opt}$, for each simulation setting. Here, $k_{opt}$ denotes the value of the deterministic $k$ that minimized the estimated error rate of the $k$-nn classifier. The Monte Carlo sampling variability of the estimated error rates can be measured by

$$\text{s.e.}(\widehat{\text{Err}}) = \sqrt{\widehat{\text{Err}}(1 - \widehat{\text{Err}})/500}.$$

It is seen that the empirical choice $\tilde{k}_{opt}$ works particularly well. The error rates of the $k$-nn classifier with $\tilde{k}_{opt}$ are not far from the error rate of the corresponding classifier with $k_{opt}$. The interval $\widehat{\text{Err}} \pm \text{s.e.}(\widehat{\text{Err}})$, where $\widehat{\text{Err}}$ is the estimated error rate of the $k$-nn classifier with $\tilde{k}_{opt}$, contains the optimal error rate achieved by the corresponding classifier with $k_{opt}$, except in the correlated case with $(\mu, \nu) = (100, 200)$. For the latter case, confidence intervals with two standard errors include the corresponding optimal value. Overall, the subsampling fraction $r = 1/3$ gave the best results. However, the error rate does not change much for different choices of $r$; the differences are not statistically significant. This suggests that $\tilde{k}_{opt}$ may not be sensitive to the choice of the resampling fraction. In the simulations we tried other populations with different mean vectors and covariance matrices. Also, we tried other training sample sizes. The lessons that we learned from the other simulation settings were basically the same as those obtained from Table 1.

Table 1 also suggests that the optimal choice $k_{opt}$ for the case $\mu \neq \nu$ tends to be smaller than the one for $\mu = \nu$. Our theory for the rate of $k_{opt}$ also was evident empirically. For example, we found that $k_{opt}$ changed from 27 to 71 when $(\mu, \nu)$ increased from $(100, 200)$ to $(400, 800)$ in the case corresponding to the bottom row of Table 1. The rate of increase in this case

TABLE 1

*Error rates of classifiers based on 500 training datasets from Poisson sample-size models with intensity $\lambda = \mu f + \nu g$, where $f$ and $g$ are densities of normal distributions as specified in the text. Here, $r$ denotes the subsampling fraction that appears in Section 2.4*

| $d$ | $(\mu, \nu)$ | $\rho$ | Bayes | $k_{opt}$ | $k$-nn with $k_{opt}$ | $k$-nn with $\tilde{k}_{opt}$ | | |
|---|---|---|---|---|---|---|---|---|
| | | | | | | $r = 1/3$ | $r = 1/2$ | $r = 2/3$ |
| 1 | (100, 100) | | 0.3072 | 103 | 0.3119 | 0.3119 | 0.3118 | 0.3120 |
| | (100, 200) | | 0.2685 | 61 | 0.2735 | 0.2759 | 0.2784 | 0.2814 |
| 2 | (100, 100) | 0 | 0.2371 | 71 | 0.2444 | 0.2445 | 0.2450 | 0.2454 |
| | | 0.5 | 0.1566 | 39 | 0.1654 | 0.1682 | 0.1708 | 0.1731 |
| | (100, 200) | 0 | 0.2125 | 45 | 0.2199 | 0.2236 | 0.2274 | 0.2310 |
| | | 0.5 | 0.1430 | 27 | 0.1514 | 0.1684 | 0.1784 | 0.1870 |

was $71/27 = 2.63$, which was roughly consistent with the theoretical value $4^{4/(2+4)} = 2.52$. To obtain similar empirical evidence in higher-dimensional feature spaces, we considered a case with $d = 16$. We simulated 500 training datasets from Poisson sample-size models, with $f$ and $g$ being the densities of $N_{16}((0.25, \ldots, 0.25), I_{16})$ and $N_{16}((-0.25, \ldots, -0.25), I_{16})$, respectively, when $(\mu, \nu) = (100, 200)$ and $(10000, 20000)$. The relative increase of $k_{\text{opt}}$ in this case was $61/25 = 2.44$, which is not far from its theoretical value $100^{4/(16+4)} = 2.51$.

**4. Proof of Theorem 1.** Let $\mathcal{S}^{\varepsilon}$ denote the set of points in $\mathcal{R}$ that are distant no further than $\varepsilon > 0$ from $\mathcal{S}$. Write $\mathcal{R} \setminus \mathcal{S}^{\varepsilon}$ for the set of points in $\mathcal{R}$ that are not in $\mathcal{S}^{\varepsilon}$. Using Markov's inequality, it can be shown that, for each fixed $C, \varepsilon > 0$, we have, as $\nu \to \infty$,

$$(4.1) \quad P^{\text{Pois}}(z \text{ classified by } k\text{-nn rule as type } X) = I\{\psi(z) > \tfrac{1}{2}\} + O(\nu^{-C}),$$

uniformly in $z \in \mathcal{R} \setminus \mathcal{S}^{\varepsilon}$. By letting $\varepsilon = \varepsilon(\nu)$ converge to zero sufficiently slowly in (4.1), we ensure that that result remains true for decreasingly small $\varepsilon$. We need (4.1) only when $C = 1$, and for $\varepsilon(\nu)$ decreasing sufficiently slowly to zero. This version of (4.1) implies that

$$(4.2) \quad \begin{aligned} &\int_{\mathcal{R} \setminus \mathcal{S}^{\varepsilon}} g(z) P^{\text{Pois}}(z \text{ classified by } k\text{-nn rule as type } X) \\ &\qquad = \int_{\mathcal{R} \setminus \mathcal{S}^{\varepsilon}} g(z) I\{\psi(z) > \tfrac{1}{2}\} \, dz + O(\nu^{-1}). \end{aligned}$$

In view of (4.1) and (4.2), properties of $f$ and $g$ away from $\mathcal{S}^{\varepsilon}$ do not affect the size of regret up to any polynomial order. Hence, there is no loss of generality in working with distributions for which $f$ and $g$ have two continuous derivatives on $\mathcal{R}^{\varepsilon}$, rather than simply on $\mathcal{S}^{\varepsilon}$. This simplifies notation, and so we shall make the assumption below without further comment.

Given $z \in \mathcal{R}$, let $Z_{(1)}, Z_{(2)}, \ldots$ denote the point locations of the marked point process $\mathcal{P}$, ordered such that $\|z - Z_{(1)}\| \le \|z - Z_{(2)}\| \le \cdots$; let $z_{(1)}$, $z_{(2)}, \ldots$ represent particular values of $Z_{(1)}, Z_{(2)}, \ldots$, respectively; and put $\vec{z} = (z_{(1)}, \ldots, z_{(k)})$ and $\vec{Z} = (Z_{(1)}, \ldots, Z_{(k)})$. Denote by $\Pi^{\text{Pois}}(\vec{z}, k)$ the probability, conditional on $\vec{Z} = \vec{z}$, that among the points $z_{(1)}, \ldots, z_{(k)}$ there are at least $\tfrac{1}{2}k$ points with marks $X$. We may write

$$\Pi^{\text{Pois}}(\vec{z}, k) = P\left(\sum_{i=1}^{k} J_i \ge \tfrac{1}{2}k\right),$$

where $J_1, \ldots, J_k$ are independent zero-one variables,

$$(4.3) \quad P(J_i = 1) = q_i \equiv \psi(z_{(i)}), \qquad \psi = \frac{\mu f}{\mu f + \nu g}.$$

To aid interpretation of (4.3), note that, since we are here conditioning on $\vec{Z} = \vec{z}$, $P(J_i = 1) = P(J_i = 1 \mid Z_{(i)} = z_{(i)})$.

Note that, uniformly in $1 \leq i \leq k \in [k_1(\nu), k_2(\nu)]$,

$$
\begin{aligned}
E\{\psi(Z_{(i)})\} = {} & \psi(z) + \sum_{j=1}^{d} E(Z_{(i)} - z)^{(j)} \psi_j(z) \\
& + \frac{1}{2} \sum_{j_1=1}^{d} \sum_{j_2=1}^{d} E\{(Z_{(i)} - z)^{(j_1)}(Z_{(i)} - z)^{(j_2)}\} \psi_{j_1 j_2}(z) \\
& + o(E\|Z_{(k)} - z\|^2),
\end{aligned}
\tag{4.4}
$$

where $(Z_{(i)} - z)^{(j)}$ denotes the $j$th component of $Z_{(i)} - z$, $\psi_j(z) = (\partial/\partial z^{(j)})\psi(z)$ and $\psi_{j_1 j_2}(z) = (\partial^2/\partial z^{(j_1)} \partial z^{(j_2)})\psi(z)$. To obtain (4.4), we have used (2.4), which implies that, for sufficiently small $\varepsilon > 0$, $f$ and $g$ have two continuous derivatives on $\mathcal{R}^\varepsilon$, the latter denoting the set of all points in $\mathbb{R}^d$ that are distant no further than $\varepsilon$ from some point in $\mathcal{R}$. It follows from this result that, under (2.4), the probability that $Z_i = Z_i(z) \in \mathcal{R}^\varepsilon$, for all $1 \leq i \leq k_2(\nu)$ and all $z \in \mathcal{R}$, equals $1 - O(\nu^{-C})$ for all $C > 0$. This implies the Taylor expansion of $\psi(Z_{(i)})$ that leads to (4.4), and, in combination with the moment condition in (2.4), ensures the correctness of the remainder term in (4.4).

Under the conditions of Theorem 1, $E\|Z_{(k)} - z\|^2 = O\{(k/\nu)^{2/d}\}$, and so (4.4) implies that

$$
\begin{aligned}
\sum_{i=1}^{k} & \{E\psi(Z_{(i)}) - \psi(z)\} \\
& = \sum_{i=1}^{k} \sum_{j=1}^{d} E(Z_{(i)} - z)^{(j)} \psi_j(z) \\
& \quad + \frac{1}{2} \sum_{i=1}^{k} \sum_{j_1=1}^{d} \sum_{j_2=1}^{d} E\{(Z_{(i)} - z)(Z_{(i)} - z)^{\mathrm{T}}\}_{j_1 j_2} \psi_{j_1 j_2}(z) \\
& \quad + o\{k(k/\nu)^{2/d}\}.
\end{aligned}
\tag{4.5}
$$

Since $\mathcal{R}$ is compact and the remainder in (4.5) is of the stated order for each $z \in \mathcal{R}$, then the remainder is of that order uniformly in $z$.

Writing $\tau = (\mu/\nu)f + g$ and $\kappa(u, z) = \int_{v:\|v\| \leq \|u\|} \tau(z + v)\,dv$, we see that the density of $Z_{(i)} - z$ at $u$ is

$$
f_i(u, z) = \nu\tau(z + u)\frac{\{\nu\kappa(u, z)\}^{i-1}}{(i-1)!} e^{-\nu\kappa(u, z)}.
$$

Therefore,

$$\sum_{i=1}^{k} E(Z_{(i)} - z) = \nu \int u\tau(z+u)P\{W(u,z) \le k-1\}\, du, \tag{4.6}$$

$$\sum_{i=1}^{k} E\{(Z_{(i)} - z)(Z_{(i)} - z)^{\mathrm{T}}\} = \nu \int uu^{\mathrm{T}}\tau(z+u)P\{W(u,z) \le k-1\}\, du, \tag{4.7}$$

where the random variable $W(u,z)$ is Poisson-distributed with mean $\nu\kappa(u,z)$, and the integrals are over $\mathbb{R}^d$.

In (4.6) and (4.7) we shall make the change of variable

$$u = \left\{\frac{k}{\nu a_d \tau(z)}\right\}^{1/d} v. \tag{4.8}$$

If $\varepsilon_1 > 0$ is chosen so small that $\nu^{-3\varepsilon_1}\{\nu/k_2(\nu)\}^{2/d} \to \infty$, then, with $v$ defined by (4.8), and $t_j = \nu^{-j\varepsilon_1}\{\nu/k_2(\nu)\}^{2/d}$, we have, for all sufficiently large $\nu$, and for all $\|v\| > t_1$ and all $k \in [k_1(\nu), k_2(\nu)]$,

$$\frac{k-1-E\{W(u,z)\}}{\{EW(u,z)\}^{1/2}} \le \frac{k-1-kt_2}{(kt_2)^{1/2}} \le -\frac{1}{2}(kt_2)^{1/2} = -\frac{1}{2}(k\nu^{\varepsilon_1}t_3)^{1/2}.$$

It follows that, for all sufficiently large $\nu$, and for all $\|v\| \ge t_1$, all $k \in [k_1(\nu), k_2(\nu)]$ and each $C > 0$,

$$\begin{aligned}
P\{W(u,z) \le k-1\} &\le P\left\{-\frac{W(u,z)-EW(u,z)}{\{EW(u,z)\}^{1/2}} \ge \frac{1}{2}(k\nu^{\varepsilon_1}t_3)^{1/2}\right\} \\
&\le (4/\nu^{\varepsilon_1})^C E\left[\left|\frac{W(u,z)-EW(u,z)}{\{EW(u,z)\}^{1/2}}\right|^{2C}\right] \\
&\le C_1\nu^{-C\varepsilon_1},
\end{aligned} \tag{4.9}$$

where $C_1 > 0$ depends only on $C$. Here we have used the fact that $W(u,z)$ is Poisson-distributed with a mean that is bounded below by 1 for $\|v\| \ge t_1$ and large $\nu$.

Combining (4.5), (4.6), (4.7) and (4.9), and noting that the distribution of $W(u,z)$ is symmetric in $u$, we deduce that

$$\begin{aligned}
\sum_{i=1}^{k} &\{E\psi(Z_{(i)}) - \psi(z)\} \\
&= \nu \int_{u:\|v\| \le t_1} \dot{\psi}(z)^{\mathrm{T}} u\{\tau(z+u) - \tau(z)\}P\{W(u,z) \le k-1\}\, du \\
&\quad + \tfrac{1}{2}\nu \int_{u:\|v\| \le t_1} u^{\mathrm{T}}\ddot{\psi}(z)u\tau(z+u)P\{W(u,z) \le k-1\}\, du \\
&\quad + o\{k(k/\nu)^{2/d}\},
\end{aligned} \tag{4.10}$$

uniformly in $z \in \mathcal{R}$.

Writing $\dot{\tau} = (\tau_1, \ldots, \tau_d)^{\mathrm{T}}$, defining $\dot{\psi}$ and $\dot{\tau}$ analogously, defining $\ddot{\psi} = (\psi_{ij})$, a $d \times d$ matrix, and taking $\mathcal{T}$ to be the set of $v$ such that $\|v\| \le t_1$, and $\mathcal{T}'$ to be the corresponding set of $u$, given by (4.8), we deduce from (4.10) that

$$
\begin{aligned}
\sum_{i=1}^{k} & \{E\psi(Z_{(i)}) - \psi(z)\} \\
&= \nu \int_{\mathcal{T}'} \{\dot{\psi}(z)^{\mathrm{T}} u u^{\mathrm{T}} \dot{\tau}(z) + \tfrac{1}{2} u^{\mathrm{T}} \ddot{\psi}(z) u \tau(z)\} \\
&\qquad \times P\{W(u,z) \le k-1\}\, du + o\{k(k/\nu)^{2/d}\} \\
&= \nu \int_{\mathcal{T}'} \left[ \sum_{j=1}^{d} (u^{(j)})^2 \{\psi_j(z)\tau_j(z) + \tfrac{1}{2}\psi_{jj}(z)\tau(z)\} \right] \\
&\qquad \times P\{W(u,z) \le k-1\}\, du + o\{k(k/\nu)^{2/d}\} \\
&= \{k/a_d\tau(z)\}\{k/\nu a_d\tau(z)\}^{2/d} \\
&\qquad \times \int_{\mathcal{T}} \left[ \sum_{j=1}^{d} (v^{(j)})^2 \{\psi_j(z)\tau_j(z) + \tfrac{1}{2}\psi_{jj}(z)\tau(z)\} \right] \\
&\qquad \times P\{W(u,z) \le k-1\}\, dv + o\{k(k/\nu)^{2/d}\},
\end{aligned}
$$

(4.11)

uniformly in $z \in \mathcal{R}$.

To control the value of $P\{W(u,z) \le k-1\}$ in (4.11), we shall use a normal approximation to the distribution of a Poisson random variable with large mean, and a crude bound to that distribution when the mean is small. Specifically, let $Z_\zeta$ have a Poisson distribution with mean $\zeta$. Then, for each $C > 0$, there exists a constant $C_1 = C_1(C) > 0$ such that, whenever $\zeta \ge 0$,

(4.12)  (a) for $\zeta \ge 1$, $\sup_{-\infty < x < \infty} (1 + |x|)^C |P(Z_\zeta \le \zeta + \zeta^{1/2}x) - \Phi(x)| \le C_1 \zeta^{-1/2}$, and (b) for $0 \le \zeta \le 1$, $\sup_{x>0} (1 + |x|)^C P(Z_\zeta > x) \le C_1 \zeta$.

Since $\kappa(u,z) = a_d\tau(z)\|u\|^d \{1 + O(\|u\|^2)\}$ as $\|u\| \to 0$, uniformly in $z \in \mathcal{R}$, then if $u$ is given by (4.8), $\nu\kappa(u,z) = k\|v\|^d [1 + O\{(k/\nu)^{2/d}\|v\|^2\}]$, uniformly in $z \in \mathcal{R}$. It follows that

$$
\begin{aligned}
\frac{k-1-EW(u,z)}{\{EW(u,z)\}^{1/2}} &= \frac{k^{1/2}(1 - \|v\|^d)}{\|v\|^{d/2}} [1 + O\{(k/\nu)^{2/d}\|v\|^2\}] \\
&\qquad + O\{k^{1/2}(k/\nu)^{2/d}\|v\|^{2+(d/2)} + k^{-1/2}\|v\|^{-d/2}\}.
\end{aligned}
$$

Noting that

(4.13)  $P\{W(u,z) \le k-1\} = P\left[ \dfrac{W(u,z) - EW(u,z)}{\{\operatorname{var} W(u,z)\}^{1/2}} \le \dfrac{k-1-EW(u,z)}{\{EW(u,z)\}^{1/2}} \right],$

using (4.12)(a) to produce an approximation to the right-hand side of (4.13) when $k^{-1/d} \leq \|v\| \leq t_1$, and using (4.12)(b) for the same purpose when $\|v\| \leq k^{-1/d}$, we deduce from (4.11) that

$$(4.14) \qquad \sum_{i=1}^{k} \{E\psi(Z_{(i)}) - \psi(z)\} = k(k/\nu)^{2/d}\alpha_1(z) + o\{k(k/\nu)^{2/d}\},$$

uniformly in $z \in \mathcal{R}$, where

$$
\begin{aligned}
(4.15) \qquad \alpha_1(z) &\equiv \{a_d\tau(z)\}^{-1-(2/d)} \\
&\quad \times \int_{\|v\| \leq 1} \left[ \sum_{j=1}^{d} (v^{(j)})^2 \{\psi_j(z)\tau_j(z) + \tfrac{1}{2}\psi_{jj}(z)\tau(z)\} \right] dv \\
&= \{a_d\tau(z)\}^{-1-(2/d)} d^{-1} \left[ \sum_{j=1}^{d} \{\psi_j(z)\tau_j(z) + \tfrac{1}{2}\psi_{jj}(z)\tau(z)\} \right] \\
&\quad \times \int_{\|v\| \leq 1} \left\{ \sum_{j=1}^{d} (v^{(j)})^2 \right\} dv,
\end{aligned}
$$

the latter being identical to $\alpha(z)$, defined at (2.6), except that there, $\psi$ and $\tau$ are replaced by their respective limits, $\rho$ and $\lambda$.

In our proofs throughout Section 4.1, it is convenient to work not with $\mathcal{S}$ but with the locus $\mathcal{S}_\nu$ of points $z_0$ such that

$$(4.16) \qquad \frac{\mu f(z_0)}{\mu f(z_0) + \nu g(z_0)} = \frac{1}{2}.$$

[In this notation $\mathcal{S}_\infty = \lim_{\nu \to \infty} \mathcal{S}_\nu$ is the set of $z_0$ such that $\rho(z_0) = \frac{1}{2}$.] We shall suppress the subscript on $\mathcal{S}_\nu$, however, instead showing at the end of the proof [see the argument below (4.23)] that the transition from $\mathcal{S} = \mathcal{S}_\nu$ to $\mathcal{S}_\infty$ is elementary.

We wish to develop an approximation to

$$
\begin{aligned}
(4.17) \qquad K_\nu(g) &\equiv \int_{\mathcal{S}^\varepsilon} g(z) P^{\text{Pois}}(z \text{ classified by } k\text{-nn rule as type } X)\, dz \\
&= \int_{\mathcal{S}^\varepsilon} g(z) E\{\Pi^{\text{Pois}}(\vec{Z}, k)\}\, dz.
\end{aligned}
$$

If we reinterpret $J_1, \ldots, J_k$ as random variables with distributions depending on $\vec{Z}$, independent conditional on $\vec{Z}$, and satisfying $P(J_i = 1 \mid \vec{Z}) = \psi(Z_{(i)})$, then

$$(4.18) \qquad E\{\Pi^{\text{Pois}}(\vec{Z}, k)\} = P\left( \sum_{i=1}^{k} J_i \geq \tfrac{1}{2}k \right).$$

Let $\mathcal{T}_{z_0}$ be the infinite line perpendicular to $\mathcal{S}$ at $z_0$, and let $u$ denote a point on $\mathcal{T}_{z_0}$. Now, $\mathcal{T}_{z_0}$ has two "halves," one in the direction where $\psi(u)$ immediately increases above $\frac{1}{2}$ as $u$ is moved away from $z_0$, and the other where $\psi(u)$ immediately decreases below $\frac{1}{2}$. Call these $\mathcal{T}_{z_0+}$ and $\mathcal{T}_{z_0-}$, respectively. Note that $\mathcal{T}_{z_0} = \{z_0 + t\dot{\psi}(z_0) : -\infty < t < \infty\}$ and $\mathcal{T}_{z_0+} = \{z_0 + t\dot{\psi}(z_0) : 0 < t < \infty\}$.

Put $\mu_k(z) = \sum_{i \le k} E(J_i)$, $\sigma_k(z)^2 = \mathrm{var}(\sum_{i \le k} J_i)$, $W_k(z) = \{\sum_{i \le k} J_i - \mu_k(z)\}/\sigma_k(z)$ and $\chi(z) = I\{\psi(z) \le \frac{1}{2}\}$. Assume $\varepsilon \downarrow 0$ and $k_1(\nu)^{1/2}\varepsilon \to \infty$ as $\nu \to \infty$. Then,

$$(4.19) \qquad \sigma_k(z)^2 = \tfrac{1}{4}k\{1 + o(1)\}, \qquad \text{uniformly in } z \in \mathcal{S}^\varepsilon,$$

as $\nu \to \infty$. By (4.17) and (4.18),

$$K'_\nu(g) \equiv K_\nu(g) - \int_{\mathcal{S}^\varepsilon} g(z)(1-\chi)(z)\,dz$$

$$= \int_{\mathcal{S}^\varepsilon} g(z)[P\{W_k(z) > w_k(z)\} - (1-\chi)(z)]\,dz,$$

$$K'_\nu(f) \equiv K_\nu(f) - \int_{\mathcal{S}^\varepsilon} f(z)\chi(z)\,dz$$

$$= \int_{\mathcal{S}^\varepsilon} f(z)[P\{W_k(z) \le w_k(z)\} - \chi(z)]\,dz,$$

where

$$w_k(z) = -\frac{1}{\sigma_k(z)} \sum_{i=1}^{k} \Big\{ E\psi(Z_{(i)}) - \frac{1}{2} \Big\}.$$

Hence,

$$(4.20) \qquad \begin{aligned} K''_\nu &\equiv \frac{\mu}{\mu+\nu} K'_\nu(f) + \frac{\nu}{\mu+\nu} K'_\nu(g) \\ &= \int_{\mathcal{S}^\varepsilon} \left\{ \frac{\mu f(z) - \nu g(z)}{\mu+\nu} \right\} [P\{W_k(z) \le w_k(z)\} - \chi(z)]\,dz. \end{aligned}$$

In view of (4.19), a standard application of the nonuniform version of the Berry–Esseen theorem to the sum of independent random variables represented by $W_k(z)$ implies that, for each $C > 0$,

$$\sup_{z \in \mathcal{S}^\varepsilon} \sup_{-\infty < w < \infty} (1 + |w|)^C |P\{W_k(z) \le w\} - \Phi(w)| = O(k^{-1/2}).$$

Hence, (4.20) entails, for all $C > 0$,

$$(4.21) \qquad \begin{aligned} K''_\nu &= \int_{\mathcal{S}^\varepsilon} \left\{ \frac{\mu f(z) - \nu g(z)}{\mu+\nu} \right\} [\Phi\{w_k(z)\} - \chi(z)]\,dz \\ &\quad + O\left[ k^{-1/2} \int_{\mathcal{S}^\varepsilon} \left| \frac{\mu f(z) - \nu g(z)}{\mu+\nu} \right| \{1 + |w_k(z)|\}^{-C}\,dz \right]. \end{aligned}$$

Using (4.14), (4.15) and (4.19), it can be shown that, if we take $z = z_0 + k^{-1/2}u$, with $z_0 \in \mathcal{S}$ and $u$ given by $z_0 + k^{-1/2}u \in \mathcal{T}_{z_0}$, then

$$-w_k(z) = \{1 + o(1)\}2k^{-1/2}k[\psi(z) - \tfrac{1}{2} + (k/\nu)^{2/d}\alpha_1(z) + o\{(k/\nu)^{2/d}\}]$$

$$= \{1 + o(1)\}2\left[\sum_{j=1}^{d} u^{(j)}\psi_j(z_0) + k^{1/2}(k/\nu)^{2/d}\alpha_1(z_0)\right.$$

$$\left. \times o\{\|u\| + k^{1/2}(k/\nu)^{2/d}\}\right],$$

uniformly in $z \in \mathcal{S}^\varepsilon$. Hence, writing $\mathcal{U}_{z_0} = \mathcal{T}_{z_0} - z_0$, $\mathcal{U}_{z_0\pm} = \mathcal{T}_{z_0\pm} - z_0$, and $c_k = k^{1/2}(k/\nu)^{2/d}$, we obtain from (4.21)

$$
\begin{aligned}
kK_\nu'' &= \int_{\mathcal{S}} \int_{\mathcal{U}_{z_0}} \{p\dot{f}(z_0) - (1-p)\dot{g}(z_0)\}^{\mathrm{T}} \\
&\qquad\qquad \times u(\Phi[-2\{\dot{\psi}(z_0)^{\mathrm{T}}u + c_k\alpha_1(z_0)\}] - I(u \in \mathcal{U}_{z_0-}))\,du\,dz_0 \\
&\quad + o(1 + c_k^2) \\
\text{(4.22)}\qquad &= \int_{\mathcal{S}} \int_{-\infty}^{\infty} \{p\dot{f}(z_0) - (1-p)\dot{g}(z_0)\}^{\mathrm{T}}\dot{\psi}(z_0)\|\dot{\psi}(z_0)\|^{-1} \\
&\qquad\qquad \times t(\Phi[-2\{\|\dot{\psi}(z_0)\|t + c_k\alpha_1(z_0)\}] - I(t < 0))\,dt\,dz_0 \\
&\quad + o(1 + c_k^2) \\
&= C_1(\mathcal{S}) + C_2(\mathcal{S})c_k^2 + o(1 + c_k^2),
\end{aligned}
$$

where, to obtain the second identity, we take $u = t\dot{\psi}(z_0)/\|\dot{\psi}(z_0)\|$. In (4.22), $C_1(\mathcal{S})$ and $C_2(\mathcal{S})$ have the definitions at (2.7), except that here $\mathcal{S}$ is interpreted as the set of points $z_0$ for which (4.16) holds.

Combining (4.2) and (4.22), we deduce that

$$\text{(4.23)}\quad \text{risk}_{k\text{-nn}}^{\text{Pois}} - \text{risk}_{\text{Bayes}}^{\text{Pois}} = C_1(\mathcal{S})k^{-1} + C_2(\mathcal{S})(k/\nu)^{4/d} + o\{k^{-1} + (k/\nu)^{4/d}\}.$$

Under the conditions assumed in Theorem 1, $\mu/(\mu+\nu) \to p$ as $\nu \to \infty$, from which it follows that $C_1(\mathcal{S})$ and $C_2(\mathcal{S})$ converge to the values they would take if we were to define $\mathcal{S}$ as the set of points $z_0$ for which, instead of (4.16), $pf(z_0)/\{pf(z_0) + (1-p)g(z_0)\} = \tfrac{1}{2}$. This is the definition used for $\mathcal{S}$ at (2.7). Note too that $\psi \to \rho$ and $\tau \to \lambda$ as $\nu \to \infty$, and that these limits arise in a very simple way. For example, $\tau = (\mu/\nu)f + g$ converges to $\lambda = p(1-p)^{-1}f + g$ since $\mu/\nu \to p(1-p)^{-1}$; the functions $f$ and $g$ remain fixed. Since $C_1(\mathcal{S})$ and $C_2(\mathcal{S})$ converge to their values at (2.7), then Theorem 1 follows from (4.23).

**Acknowledgments.** We are grateful to the reviewers for helpful comments.

P. HALL
DEPARTMENT OF MATHEMATICS AND STATISTICS
UNIVERSITY OF MELBOURNE
PARKVILLE, VIC 3010
AUSTRALIA

B. U. PARK
DEPARTMENT OF STATISTICS
SEOUL NATIONAL UNIVERSITY
SEOUL 151–747
KOREA
E-MAIL: bupark2000@gmail.com

R. J. SAMWORTH
STATISTICAL LABORATORY
CENTRE FOR MATHEMATICAL SCIENCES
UNIVERSITY OF CAMBRIDGE
WILBERFORCE ROAD
CAMBRIDGE, CB3 0WB
UNITED KINGDOM